# P-values for classification

**Lutz Dümbgen**[*]

*University of Bern*
*Institute for Mathematical Statistics and Actuarial Science*
*e-mail:* duembgen@stat.unibe.ch

**Bernd-Wolfgang Igl**[†]

*University at Lübeck*
*Institute of Medical Biometry and Statistics*
*e-mail:* bernd.igl@gmx.de

**Axel Munk**[‡]

*Georgia Augusta University Göttingen*
*Institute for Mathematical Stochastics*
*e-mail:* munk@math.uni-goettingen.de

**Abstract:** Let $(X, Y)$ be a random variable consisting of an observed feature vector $X \in \mathcal{X}$ and an unobserved class label $Y \in \{1, 2, \ldots, L\}$ with unknown joint distribution. In addition, let $\mathcal{D}$ be a training data set consisting of $n$ completely observed independent copies of $(X, Y)$. Usual classification procedures provide point predictors (classifiers) $\widehat{Y}(X, \mathcal{D})$ of $Y$ or estimate the conditional distribution of $Y$ given $X$. In order to quantify the certainty of classifying $X$ we propose to construct for each $\theta = 1, 2, \ldots, L$ a p-value $\pi_\theta(X, \mathcal{D})$ for the null hypothesis that $Y = \theta$, treating $Y$ temporarily as a fixed parameter. In other words, the point predictor $\widehat{Y}(X, \mathcal{D})$ is replaced with a prediction region for $Y$ with a certain confidence. We argue that (i) this approach is advantageous over traditional approaches and (ii) any reasonable classifier can be modified to yield nonparametric p-values. We discuss issues such as optimality, single use and multiple use validity, as well as computational and graphical aspects.

Received June 2008.

## 1. Introduction

Let $(X, Y)$ be a random variable consisting of a feature vector $X \in \mathcal{X}$ and a class label $Y \in \Theta := \{1, \ldots, L\}$ with $L \geq 2$ possible values. The joint distribution of $X$ and $Y$ is determined by the prior probabilities $w_\theta := \mathbb{P}(Y = \theta)$ and the conditional distributions $P_\theta := \mathcal{L}(X \,|\, Y = \theta)$ for all $\theta \in \Theta$. Classifying such an observation $(X, Y)$ means that only $X$ is observed, while $Y$ has to be predicted somehow. There is a vast literature on classification, and we refer to McLachlan [7], Ripley [10] or Fraley and Raftery [4] for an introduction and further references.

[*]Work supported by Swiss National Science Foundation (SNF)
[†]Work supported by German Ministry of Education and Research (BMBF)
[‡]Work supported by German Science Foundation (DFG)

Let us assume for the moment that the joint distribution of $X$ and $Y$ is known, so that training data are not needed yet. In the simplest case, one chooses a classifier $\widehat{Y} : \mathcal{X} \to \Theta$, i.e. a point predictor of $Y$. A possible extension is to consider $\widehat{Y} : \mathcal{X} \to \{0\} \cup \Theta$, where $\widehat{Y}(X) = 0$ means that no class is viewed as plausible. A Bayesian approach would be to calculate the posterior distribution of $Y$ given $X$, i.e. the posterior weights $w_\theta(X) := \mathbb{P}(Y = \theta \,|\, X)$. In fact, a classifier $\widehat{Y}^*$ satisfying

$$\widehat{Y}^*(X) \ \in \ \arg\max_{\theta \in \Theta} w_\theta(X)$$

is well-known [7, Chapter 1] to minimize the risk

$$R(\widehat{Y}) \ := \ \mathbb{P}(\widehat{Y}(X) \ne Y).$$

An obvious advantage of using the posterior distribution instead of the simple classifier $\widehat{Y}^*$ (or $\widehat{Y}$) is additional information about confidence. That means, for instance, the possibility of computing the conditional risk $\mathbb{P}(\widehat{Y}^*(X) \ne Y \,|\, X) = 1 - \max_\theta w_\theta(X)$. However, this depends very sensitively on the prior weights $w_\theta$. Small changes in the latter may result in drastic changes of the posterior weights $w_\theta(X)$. Moreover, if some classes $\theta$ have very small prior weight, the classifier $\widehat{Y}^*$ tends to ignore these, i.e. the class-dependent risk $\mathbb{P}(\widehat{Y}^*(X) \ne Y \,|\, Y = \theta)$ may be rather large for some classes $\theta$. For instance, in medical applications each class may correspond to a certain disease status while the feature vector contains information about patients, including certain symptoms. Here it would be unacceptable to classify each person as being healthy, just because the diseases in question are extremely rare. Note also that some study designs (e.g. case-control studies) allow for the estimation of the $P_\theta$ but not the $w_\theta$. Moreover, there are applications in which the $w_\theta$ change over time while it is still plausible to assume fixed conditional distributions $P_\theta$.

Another drawback of the posterior probabilities $w_\theta(X)$ is the following: Suppose that the prior weights $w_\theta$ are all identical and that for some subset $\Theta_o$ of $\Theta$ with at least two elements the conditional distributions $P_\theta$, $\theta \in \Theta_o$, are very similar. Then the posterior distribution of $Y$ given $X$ divides the mass corresponding to $\Theta_o$ essentially uniformly among its elements. Even if the point $X$ is right in the 'center' of the distributions $P_\theta$, $\theta \in \Theta_o$, so that each class in $\Theta_o$ is perfectly plausible, the posterior weights are not greater than $1/\#\Theta_o$. If $w_\theta(X)$ is viewed merely as a measure of plausibility of class $\theta$, there is no compelling reason why these measures should add to one.

To treat all classes impartially, we propose to compute for each class $\theta \in \Theta$ a p-value $\pi_\theta(X)$ of the null hypothesis that $Y = \theta$. (In this formulation we treat $Y$ temporarily as an unknown fixed parameter.) That means, $\pi_\theta : \mathcal{X} \to [0, 1]$ satisfies

$$\mathbb{P}\bigl(\pi_\theta(X) \le \alpha \,\big|\, Y = \theta\bigr) \ \le \ \alpha \quad \text{for all } \alpha \in (0, 1). \tag{1.1}$$

Given such p-values $\pi_\theta$, the set

$$\widehat{\mathcal{Y}}_\alpha(X) := \{\theta \in \Theta : \pi_\theta(X) > \alpha\}$$

is a $(1-\alpha)$–prediction region for $Y$, i.e.

$$\mathbb{P}(Y \in \widehat{\mathcal{Y}}_\alpha(X) \,|\, Y = \theta) \geq 1 - \alpha \quad \text{for arbitrary } \theta \in \Theta, \alpha \in (0,1).$$

If $\widehat{\mathcal{Y}}_\alpha(X)$ happens to be a singleton, we have classified $X$ uniquely with given confidence $1-\alpha$. In case of $2 \leq \#\widehat{\mathcal{Y}}_\alpha(X) < L$ we can at least exclude some classes with a certain confidence.

So far the classification problem corresponds to a simple statistical model with finite parameter space $\Theta$. A distinguishing feature of classification problems is that the joint distribution of $(X,Y)$ is typically unknown and has to be estimated from a set $\mathcal{D}$ consisting of completely observed training observations $(X_1, Y_1), (X_2, Y_2), \ldots, (X_n, Y_n)$. Let us assume for the moment that all $n+1$ observations, i.e. the $n$ training observations $(X_i, Y_i)$ and the current observation $(X, Y)$, are independent and identically distributed. Now one has to consider classifiers $\widehat{Y}(X, \mathcal{D})$ and p-values $\pi_\theta(X, \mathcal{D})$ depending on the current feature vector $X$ as well as on the training data $\mathcal{D}$. In this situation one can think of two possible extensions of (1.1): For any $\theta \in \Theta$ and $\alpha \in (0,1)$,

$$\mathbb{P}(\pi_\theta(X, \mathcal{D}) \leq \alpha \,|\, Y = \theta) \leq \alpha, \tag{1.2}$$

$$\mathbb{P}(\pi_\theta(X, \mathcal{D}) \leq \alpha \,|\, Y = \theta, \mathcal{D}) \leq \alpha + o_p(1) \quad \text{as } n \to \infty. \tag{1.3}$$

It will turn out that Condition (1.2) can be guaranteed in various settings. Condition (1.3) corresponds to "multiple use" of our p-values: Suppose that we use the training data $\mathcal{D}$ to construct the p-values $\pi_\theta(\cdot, \mathcal{D})$ and classify *many* future observations $(\widetilde{X}, \widetilde{Y})$. Then the relative number of future observations with $\widetilde{Y} = b$ and $\pi_\theta(\widetilde{X}, \mathcal{D}) \leq \alpha$ is close to

$$w_b \cdot \mathbb{P}(\pi_\theta(X, \mathcal{D}) \leq \alpha \,|\, Y = b, \mathcal{D}),$$

a random quantity depending on the training data $\mathcal{D}$.

P-values as discussed here have been used in some special cases before. For instance, McLachlan's [7] "typicality indices" are just p-values $\pi_\theta(X, \mathcal{D})$ satisfying (1.2) in the special case of multivariate gaussian distributions $P_\theta$; see also Section 3. However, McLachlan's p-values are used primarily to identify observations not belonging to any of the given classes in $\Theta$. In particular, they are not designed and optimized for distinguishing between classes within $\Theta$. Also the use of receiver operating characteristic (ROC) curves in the context of logistic regression or Fisher's [3] linear discriminant analysis is related to the present concept. One purpose of this paper is to provide a solid foundation for procedures of this type.

The remainder of this paper is organized as follows: In Section 2 we return to the idealistic situation of known prior weights $w_\theta$ and distributions $P_\theta$. Here we devise p-values that are optimal in a certain sense and related to the optimal classifier mentioned previously. These p-values serve as a gold standard for p-values in realistic settings. In addition we describe briefly McLachlan's [7] typicality indices and a potential compromise between the these p-values and the optimal ones.

Section 3 is devoted to p-values involving training data. After some general remarks on cross-validation and graphical representations, we discuss McLachlan's [7] p-values in view of (1.2) and (1.3). Nonparametric p-values satisfying (1.2) without any further assumptions on the distributions $P_\theta$ are proposed in Section 3.3. These p-values are based on permutation testing, and the only practical restriction is that the group sizes $N_\theta := \#\{i : Y_i = \theta\}$ within the training data should exceed the reciprocal of the intended test level $\alpha$. We claim that any reasonable classification method can be converted to yield p-values. In particular, we introduce p-values based on a suitable variant of the nearest-neighbor method. Section 3.4 deals with asymptotic properties of various p-values as the size $n$ of $\mathcal{D}$ tends to infinity. It is shown in particular that under mild regularity conditions the nearest-neighbor p-values are asymptotically equivalent to the optimal methods of Section 2. These results are analogous to results of Stone [12, Section 8] for nearest-neighbor classifiers. In Section 3.5 the nonparametric p-values are illustrated with simulated and real data. Finally, in Section 3.6 we comment on Condition (1.3) and show that the $o_p(1)$ cannot be avoided in general.

In Section 4 we comment briefly on computational aspects of our methods. Section 5 introduces the notion of 'local identifiability' for finite mixtures, which is of independent interest. For us it is helpful to define the optimal p-values in a simple manner and it is also useful for the asymptotic considerations in Section 3.4. Proofs and technical arguments are deferred to Section 6.

Let us mention a different type of confidence procedure for classification: Suppose that $\big[a_\theta(X, \mathcal{D}), b_\theta(X, \mathcal{D})\big]$ is a confidence interval for $w_\theta(X)$. Precisely, let $a_\theta(X, \mathcal{D}) \leq w_\theta(X) \leq b_\theta(X, \mathcal{D})$ for all $\theta \in \Theta$ with probability at least $1 - \alpha$. Then

$$\check{\mathcal{Y}}(X, \mathcal{D}) := \Big\{\theta \in \Theta : b_\theta(X, \mathcal{D}) \geq \max_{\eta \in \Theta} a_\eta(X, \mathcal{D})\Big\}$$

would be a prediction region for $Y$ such that $\widehat{Y}^*(X) \subset \check{\mathcal{Y}}(X, \mathcal{D})$ with probability at least $1 - \alpha$. Note, however, that this gives no control over the probability that $Y \notin \check{\mathcal{Y}}(X, \mathcal{D})$. In fact, the latter probability could be close to 50 percent. By way of contrast, with the p-values in the present paper we can guarantee to cover $Y$ with a certain confidence, even in situations where consistent estimation of the conditional probabilities $w_\theta(X)$ is difficult or even impossible.

## 2. Optimal p-values and alternatives

Suppose that the distributions $P_1, \ldots, P_L$ have known densities $f_1, \ldots, f_L > 0$ with respect to some measure $M$ on $\mathcal{X}$. Then the marginal distribution of $X$ has density $f := \sum_{b \in \Theta} w_b f_b$ with respect to $M$, and

$$w_\theta(x) = \frac{w_\theta f_\theta(x)}{f(x)}.$$

Hence the optimal classifier $\widehat{Y}^*$ may be characterized by

$$\widehat{Y}^*(X) \in \arg\max_{\theta \in \Theta} w_\theta f_\theta(X).$$

### 2.1. Optimal p-values

Here is an analogous consideration for p-values. Let $\boldsymbol{\pi} = (\pi_\theta)_{\theta \in \Theta}$ consist of p-values $\pi_\theta$ satisfying (1.1). Given the latter constraint, our goal is to provide small p-values and small predicion regions. Hence two natural measures of risk are, for instance,

$$\mathcal{R}(\boldsymbol{\pi}) := \mathbb{E} \sum_{\theta \in \Theta} \pi_\theta(X) \quad \text{or} \quad \mathcal{R}_\alpha(\boldsymbol{\pi}) := \mathbb{E} \#\widehat{\mathcal{Y}}_\alpha(X).$$

Elementary calculations reveal that

$$\mathcal{R}(\boldsymbol{\pi}) = \int_0^1 \mathcal{R}_\alpha(\boldsymbol{\pi}) \, d\alpha \quad \text{and} \quad \mathcal{R}_\alpha(\boldsymbol{\pi}) = \sum_{\theta \in \Theta} \mathcal{R}_\alpha(\pi_\theta)$$

with

$$\mathcal{R}_\alpha(\pi_\theta) := \mathbb{P}(\pi_\theta(X) > \alpha).$$

Thus we focus on minimizing $\mathcal{R}_\alpha(\pi_\theta)$ for arbitrary fixed $\theta \in \Theta$ and $\alpha \in (0,1)$ under the constraint (1.1). Since $x \mapsto 1\{\pi_\theta(x) > \alpha\}$ may be viewed as a level–$\alpha$ test of $P_\theta$ versus $\sum_{b \in \Theta} w_b P_b$, a straightforward application of the Neyman-Pearson Lemma shows that the p-value

$$\pi_\theta^*(x) := P_\theta\{z \in \mathcal{X} : (f_\theta/f)(z) \leq (f_\theta/f)(x)\}$$

is optimal, provided that the distribution $\mathcal{L}((f_\theta/f)(X))$ is continuous. Two other representations of $\pi_\theta^*$ are given by

$$\pi_\theta^*(x) = P_\theta\{z \in \mathcal{X} : w_\theta(z) \leq w_\theta(x)\}$$
$$= P_\theta\{z \in \mathcal{X} : T_\theta^*(z) \geq T_\theta^*(x)\}$$

with $T_\theta^* := \sum_{b \neq \theta} w_{b,\theta} f_b/f_\theta$ and $w_{b,\theta} := w_b/\sum_{c \neq \theta} w_c$. The former representation shows that $\pi_\theta^*(x)$ is a non-decreasing function of $w_\theta(x)$. The latter representation shows that the prior weight $w_\theta$ itself is irrelevant for the optimal p-value $\pi_\theta^*(x)$; only the ratios $w_c/w_b$ with $b, c \neq \theta$ matter. In particular, in case of $L = 2$ classes, the optimal p-values do not depend on the prior distribution of $Y$ at all.

Here and throughout this paper we assume the likelihood ratios $T_\theta^*(X)$ to have a continuous distribution. It will be shown in Section 5 that many standard families of distributions fulfill this condition. In particular, it is satisfied in case of $\mathcal{X} = \mathbb{R}^q$ and $P_\theta = \mathcal{N}_q(\mu_\theta, \Sigma_\theta)$ with parameters $(\mu_\theta, \Sigma_\theta)$, $\Sigma_\theta$ nonsingular, not all being identical. Further examples include the multivariate $t$-family as it has been advocated by Peel and McLachlan [8] to robustify cluster and discriminant analysis. These authors also discuss maximum likelihood via the EM algorithm in this model. Without the continuity condition on $\mathcal{L}(T_\theta^*(X))$ one could still devise optimal p-values by introducing *randomized* p-values, but we refrain from such extensions.

Let us illustrate the optimal p-values in two examples involving normal distributions:

**Example 2.1. (Standard model)** Let $P_\theta = \mathcal{N}_q(\mu_\theta, \Sigma)$ with mean vectors $\mu_\theta \in \mathbb{R}^q$ and a common symmetric, nonsingular covariance matrix $\Sigma \in \mathbb{R}^{q \times q}$. Then

$$T_\theta^*(x) = \sum_{b \neq \theta} w_{b,\theta} \exp\bigl((x - \mu_{\theta,b})^\top \Sigma^{-1}(\mu_b - \mu_\theta)\bigr) \qquad (2.1)$$

with $\mu_{\theta,b} := 2^{-1}(\mu_\theta + \mu_b)$. In the special case of $L = 2$ classes, let $Z(x) := (x - \mu_{1,2})^\top \Sigma^{-1}(\mu_2 - \mu_1)/\|\mu_1 - \mu_2\|_\Sigma$ with the Mahalanobis norm $\|v\|_\Sigma := (v^\top \Sigma^{-1} v)^{1/2}$. Then elementary calculations show that

$$\pi_1^*(x) = \Phi\bigl(-Z(x) - \|\mu_1 - \mu_2\|_\Sigma/2\bigr),$$
$$\pi_2^*(x) = \Phi\bigl(+Z(x) - \|\mu_1 - \mu_2\|_\Sigma/2\bigr),$$

where $\Phi$ denotes the standard gaussian c.d.f.. In case of $\|\mu_1 - \mu_2\|_\Sigma/2 \geq \Phi^{-1}(1-\alpha)$,

$$\widehat{\mathcal{Y}}_\alpha(x) = \begin{cases} \{1\} & \text{if } Z(x) < -\|\mu_1 - \mu_2\|_\Sigma/2 + \Phi^{-1}(1-\alpha), \\ \{2\} & \text{if } Z(x) > +\|\mu_1 - \mu_2\|_\Sigma/2 - \Phi^{-1}(1-\alpha), \\ \emptyset & \text{else.} \end{cases}$$

Thus the two classes are separated well so that any observation $X$ is classified uniquely (or viewed as suspicious) with confidence $1 - \alpha$. In case of $\|\mu_1 - \mu_2\|_\Sigma/2 < \Phi^{-1}(1-\alpha)$, the feature space contains regions with unique prediction and a region in which both class labels are plausible:

$$\widehat{\mathcal{Y}}_\alpha(x) = \begin{cases} \{1\} & \text{if } Z(x) \leq +\|\mu_1 - \mu_2\|_\Sigma/2 - \Phi^{-1}(1-\alpha), \\ \{2\} & \text{if } Z(x) \geq -\|\mu_1 - \mu_2\|_\Sigma/2 + \Phi^{-1}(1-\alpha), \\ \{1, 2\} & \text{else.} \end{cases}$$

**Example 2.2.** Consider $L = 3$ classes with equal prior weights $w_\theta = 1/3$ and bivariate normal distributions $P_\theta = \mathcal{N}_2(\mu_\theta, \Sigma_\theta)$, where

$$\mu_1 = (-1, 1)^\top, \quad \mu_2 = (-1, -1)^\top, \quad \mu_3 = (2, 0)^\top$$

and

$$\Sigma_1 = \Sigma_2 = \begin{pmatrix} 1 & 1/2 \\ 1/2 & 1 \end{pmatrix}, \quad \Sigma_3 = \begin{pmatrix} 0.4 & 0 \\ 0 & 0.4 \end{pmatrix}.$$

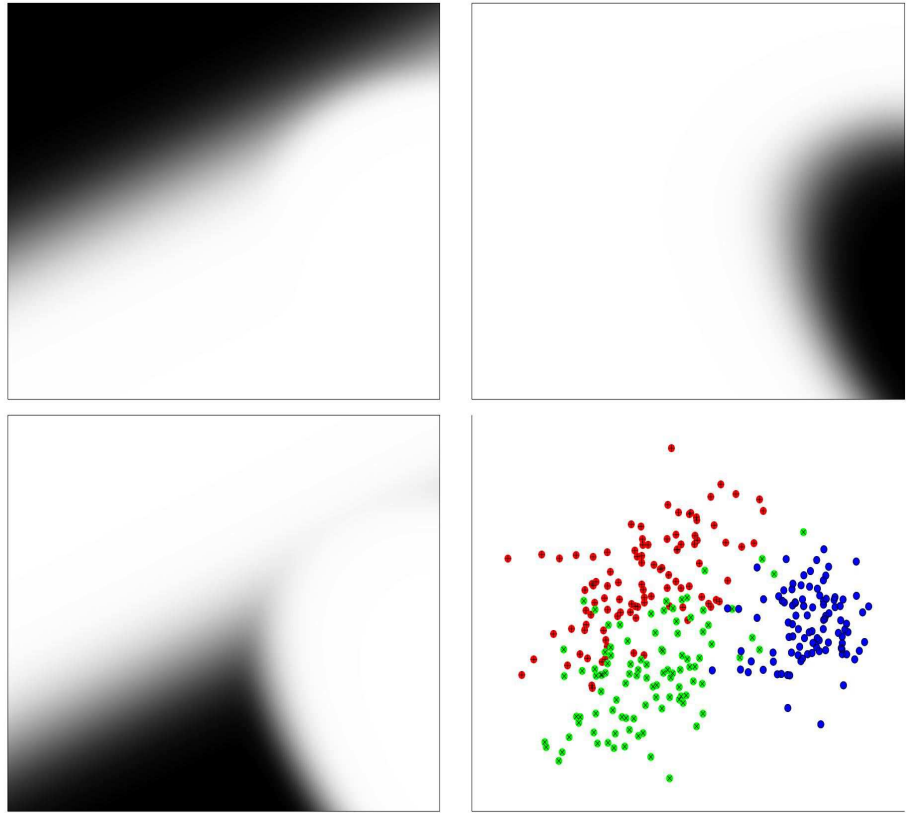

Fig 1. *P-value functions $\pi_1^*$ (top left), $\pi_2^*$ (bottom left), $\pi_3^*$ (top right) and a typical data set (bottom right) for Example 2.2.*

Figure 1 shows a typical sample from this distribution and the corresponding p-value functions $\pi_\theta^*$. The latter are on a grey scale with white corresponding to zero and black corresponding to one. The resulting predition regions $\widehat{\mathcal{Y}}_\alpha(x)$ for $\alpha = 5\%$ and $\alpha = 1\%$ are depicted in Figure 2. In the latter plots, the color of a point $x \in \mathbb{R}^2$ has the following meaning:

| Color | $\widehat{\mathcal{Y}}_\alpha(x)$ | Color | $\widehat{\mathcal{Y}}_\alpha(x)$ |
|---|---|---|---|
| black | $\emptyset$ | white | $\{1, 2, 3\}$ |
| red | $\{1\}$ | yellow | $\{1, 2\}$ |
| green | $\{2\}$ | cyan | $\{2, 3\}$ |
| dark blue | $\{3\}$ | magenta | $\{1, 3\}$ |

(The configuration $\widehat{\mathcal{Y}}_\alpha(x) = \{1,3\}$ never appeared.) Note the influence of $\alpha$: On the one hand, $\widehat{\mathcal{Y}}_{0.05}(x) = \emptyset$ for some $x \in \mathbb{R}^2$ but $\widehat{\mathcal{Y}}_{0.05}(\cdot) \neq \{1,2,3\}$ in the depicted rectangle. On the other hand, $\widehat{\mathcal{Y}}_{0.01}(x) = \{1,2,3\}$ for some $x \in \mathbb{R}^2$ while $\widehat{\mathcal{Y}}_{0.01}(\cdot) \neq \emptyset$.

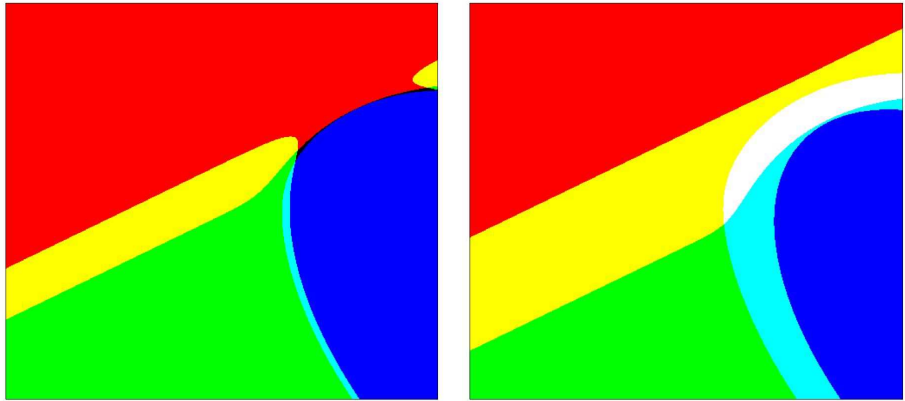

FIG 2. *Prediction regions $\widehat{\mathcal{Y}}_\alpha(x)$ for $\alpha = 5\%$ (left) and $\alpha = 1\%$ (right) in Example 2.2.*

## 2.2. Typicality indices

An alternative definition of p-values is based on the densities themselves, namely,
$$\tau_\theta(x) \ := \ P_\theta\big\{z \in \mathcal{X} : f_\theta(z) \le f_\theta(x)\big\}.$$
These typicality indices quantify to what extent a point $x$ is an outlier with respect to the single distributions $P_\theta$. These p-values $\tau_\theta$ are certainly suboptimal in terms of the risk $\mathcal{R}_\alpha(\pi_\theta)$. On the other hand, they allow for the detection of observations which belong to *none* of the classes under consideration.

**Example 2.3.** Again let $\mathcal{X} = \mathbb{R}^q$ and $P_\theta = \mathcal{N}_q(\mu_\theta, \Sigma_\theta)$. Since $f_\theta(X)$ is a strictly decreasing function of $\|X - \mu_\theta\|^2_{\Sigma_\theta}$ with conditional distribution $\chi^2_q$ given $Y = \theta$, the typicality indices may be expressed as
$$\tau_\theta(x) \ = \ 1 - F_q\big(\|x - \mu_\theta\|^2_{\Sigma_\theta}\big),$$
where $F_q$ denotes the c.d.f. of $\chi^2_q$. These p-values allow for the separation of two different classes $\theta, b \in \Theta$ only if
$$q^{-1}\|\mu_\theta - \mu_b\|^2_\Sigma$$
is sufficiently large. Thus they suffer from the curse of dimensionality and may yield much more conservative predition regions than the p-values $\pi^*_\theta$.

## 2.3. Combining the optimal p-values and typicality indices

The optimal p-values $\pi^*_\theta$ and the typicality indices $\tau_\theta$ may be viewed as extremal members of a whole family of p-values if we introduce an additional class label 0 with 'density' $f_0 \equiv 1$ and prior weight $w_0 > 0$. Then we define the compromise p-value
$$\widetilde{\pi}_\theta(x) \ := \ P_\theta\big\{z \in \mathcal{X} : (f_\theta/\widetilde{f})(z) \le (f_\theta/\widetilde{f})(x)\big\}$$

with $\widetilde{f} := \sum_{b=0}^{L} w_b f_b = f + w_0$. Note that $\widetilde{\pi}_\theta \to \tau_\theta$ pointwise as $w_0 \to \infty$, whereas $\widetilde{\pi}_\theta \to \pi_\theta^*$ as $w_0 \to 0$.

**Example 2.4.** In the setting of Example 2.1 there is another modification which is similar in spirit to Ehm et al. [1]: When defining the p-value for a particular class $\theta$ we replace the other distributions $P_b = \mathcal{N}_q(\mu_b, \Sigma)$, $b \neq \theta$, with $\widetilde{P}_b = \mathcal{N}_q(\mu_b, c\Sigma)$ for some constant $c > 1$. Thus our modified p-value becomes

$$\widetilde{\pi}_\theta(x) := P_\theta\{z \in \mathcal{X} : \widetilde{T}_\theta(z) \geq \widetilde{T}_\theta(x)\},$$

where

$$\begin{aligned}\widetilde{T}_\theta(x) &= \sum_{b=1}^{L} w_{b,\theta} \exp\bigl(\|x - \mu_\theta\|_\Sigma^2/2 - \|x - \mu_b\|_\Sigma^2/(2c)\bigr) \\ &= \sum_{b=1}^{L} w_{b,\theta} \exp\bigl((1 - c^{-1})\|x - \nu_{\theta,b}\|_\Sigma^2/2 - (c-1)^{-1}\|\mu_b - \mu_\theta\|_\Sigma^2/2\bigr)\end{aligned}$$

with $\nu_{\theta,b} := \mu_\theta - (c-1)^{-1}(\mu_b - \mu_\theta)$.

## 3. Training data

Now we return to the realistic situation of unknown distributions $P_\theta$ and p-values $\pi_\theta(X, \mathcal{D})$ with corresponding prediction regions $\widehat{\mathcal{Y}}_\alpha(X, \mathcal{D})$. From now on we consider the class labels $Y_1, Y_2, \ldots, Y_n$ as fixed while $X_1, X_2, \ldots, X_n$ and $(X, Y)$ are independent with $\mathcal{L}(X_i) = P_{Y_i}$. That way we can cover the case of i.i.d. training data (via conditioning) as well as situations with stratified training samples. In what follows let

$$\mathcal{G}_\theta := \{i \in \{1, \ldots, n\} : Y_i = \theta\} \quad \text{and} \quad N_\theta := \#\mathcal{G}_\theta.$$

We shall tacitly assume that all group sizes $N_\theta$ are strictly positive, and asymptotic statements as in (1.3) are meant as

$$n \to \infty \quad \text{and} \quad N_b/n \to w_b \quad \text{for all } b \in \Theta. \tag{3.1}$$

### 3.1. Visual assessment and estimation of separability

Before giving explicit examples of p-values, let us describe our way of visualizing the separability of different classes by means of given p-values $\pi_\theta(\cdot, \cdot)$. For that purpose we propose to compute cross-validated p-values

$$\pi_\theta(X_i, \mathcal{D}_i)$$

for $i = 1, 2, \ldots, n$ with $\mathcal{D}_i$ denoting the training data without observation $(X_i, Y_i)$. Thus each training observation $(X_i, Y_i)$ is treated temporarily as a 'future' observation to be classified with the remaining data $\mathcal{D}_i$. Then we display these cross-validated p-values graphically. This is particularly helpful for training samples of small or moderate size.

In addition to graphical displays one can compute the empirical conditional inclusion probabilities

$$\widehat{\mathcal{I}}_\alpha(b,\theta) \ := \ \#\{i \in \mathcal{G}_b : \theta \in \widehat{\mathcal{Y}}_\alpha(X_i, \mathcal{D}_i)\}/N_b$$

and the empirical pattern probabilities

$$\widehat{\mathcal{P}}_\alpha(b,S) \ := \ \#\{i \in \mathcal{G}_b : \widehat{\mathcal{Y}}_\alpha(X_i, \mathcal{D}_i) = S\}/N_b$$

for $b, \theta \in \Theta$ and $S \subset \Theta$. These numbers $\widehat{\mathcal{I}}_\alpha(b,\theta)$ and $\widehat{\mathcal{P}}_\alpha(b,S)$ can be interpreted as estimators of

$$\mathcal{I}_\alpha(b,\theta \,|\, \mathcal{D}) \ := \ \mathbb{P}\bigl(\theta \in \widehat{\mathcal{Y}}_\alpha(X, \mathcal{D}) \,\big|\, Y = b, \mathcal{D}\bigr)$$

and

$$\mathcal{P}_\alpha(b,S \,|\, \mathcal{D}) \ := \ \mathbb{P}\bigl(\widehat{\mathcal{Y}}_\alpha(X, \mathcal{D}) = S \,\big|\, Y = b, \mathcal{D}\bigr),$$

respectively; see also Section 3.4.

For large group sizes $N_b$, one can also display the empirical ROC curves

$$(0,1) \ni \alpha \ \mapsto \ 1 - \widehat{\mathcal{I}}_\alpha(b,\theta)$$

which are closely related to the usual ROC curves employed, for instance, in logistic regression or linear discriminant analysis involving $L = 2$ classes.

### 3.2. Typicality indices

For the sake of simplicity, suppose that $P_\theta = \mathcal{N}_q(\mu_\theta, \Sigma)$ with unknown mean vectors $\mu_1, \ldots, \mu_L \in \mathbb{R}^q$ and an unknown nonsingular covariance matrix $\Sigma \in \mathbb{R}^{q \times q}$. Consider the standard estimators

$$\widehat{\mu}_\theta \ := \ N_\theta^{-1} \sum_{i \in \mathcal{G}_\theta} X_i \quad \text{and} \quad \widehat{\Sigma} \ := \ (n-L)^{-1} \sum_{i=1}^n (X_i - \widehat{\mu}_{Y_i})(X_i - \widehat{\mu}_{Y_i})^\top.$$

Then the squared Mahalanobis distance

$$T_\theta(X, \mathcal{D}) \ := \ \bigl\|X - \widehat{\mu}_\theta\bigr\|_{\widehat{\Sigma}}^2$$

can be used to assess the plausibility of class $\theta$, where we assume that $n \geq L+q$. Precisely,

$$C_\theta \ := \ \frac{(n-L-q+1)}{q(n-L)(1+N_\theta^{-1})}$$

is a normalizing constant such that

$$C_\theta T_\theta(X, \mathcal{D}) \ \sim \ F_{q, n-L-q+1} \,\big|\, Y = \theta;$$

see [7]. Here $F_{k,z}$ denotes the $F$-distribution with $k$ and $z$ degrees of freedom, and we use the same symbol for the corresponding c.d.f.. Hence the typicality index

$$\tau_\theta(X, \mathcal{D}) := 1 - F_{q,n-L-q+1}(C_\theta T_\theta(X, \mathcal{D}))$$

is a p-value satisfying (1.2). Moreover, since the estimators $\widehat{\mu}_b$ and $\widehat{\Sigma}$ are consistent, one can easily verify property (1.3) as well.

**Example 3.1.** An array of ten electrochemical sensors is used for "smelling" different substances. In each case it produces raw data $\widetilde{X} \in \mathbb{R}^{10}$ consisting of the electrical resistances of these sensors. Before analyzing such data one should standardize them in order to achieve invariance with respect to the substance's concentration. One possible standardization is to replace $\widetilde{X}$ with

$$X := \left( \widetilde{X}(j) \Big/ \sum_{k=1}^{10} \widetilde{X}(k) \right)_{j=1}^{9}.$$

Thus we end up with data vectors in $\mathbb{R}^9$. For technical reasons, group sizes $N_\theta$ are typically small, and not too many future observations may be analysed. This is due to the fact that the system needs to be recalibrated regularly.

Now we consider a specific dataset with "smells" of $L = 12$ different brands of tobacco and fixed group sizes $N_\theta = 3$ for all $\theta \in \Theta$. We computed the cross-validated typicality indices $\tau_\theta(X_i, \mathcal{D}_i)$ described above. Figure 3 depicts for each training observation $(X_i, Y_i)$ the p-values $\tau_1(X_i, \mathcal{D}_i), \ldots, \tau_{12}(X_i, \mathcal{D}_i)$ as a row of twelve rectangles. The area of these is proportional to the corresponding p-value. The first three rows correspond to data from the first brand, the next three rows to the second brand, and so on. Figure 4 displays the corresponding prediction regions $\widehat{\mathcal{Y}}_\alpha(X_i, \mathcal{D}_i)$ for $\alpha = 0.01$. Within each row the elements of $\widehat{\mathcal{Y}}_\alpha(X_i, \mathcal{D}_i)$ are indicated by rectangles of full size. These pictures show classes 1 and 2 are separated well from the other eleven classes. Classes 5, 8, 9 and 12 overlap somewhat but are clearly separated from the remaining eight classes. Finally there are three pairs of classes which are essentially impossible to distinguish, at least with the present method, but which are separated well from the other ten classes. These pairs are 3-4, 6-7, and 10-11. It turned out later that brands 6 and 7 were in fact identical. Note also that all except one prediction region $\widehat{\mathcal{Y}}_\alpha(X_i, \mathcal{D}_i)$ contain the true class and at most three additional class labels.

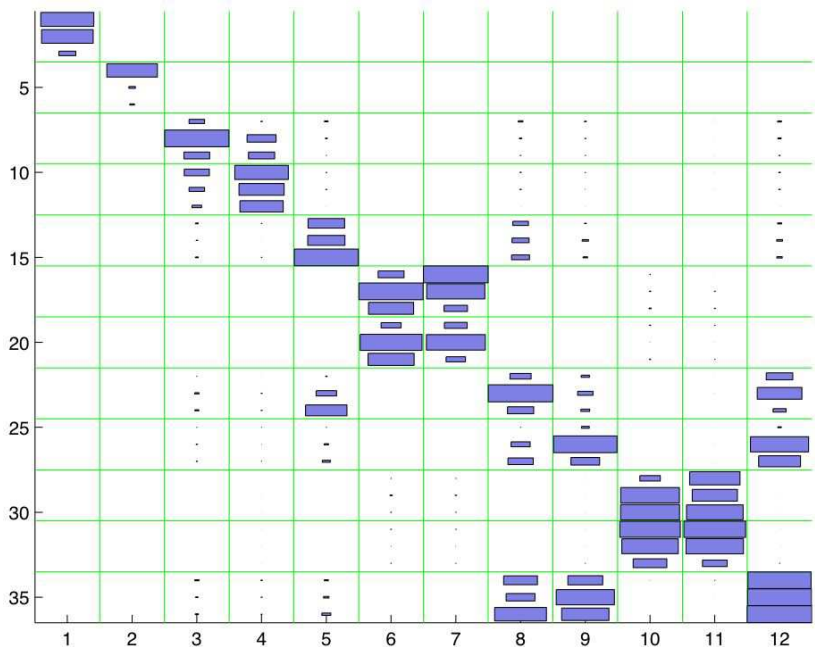

FIG 3. *Cross-validated typicality indices for tobacco "smells".*

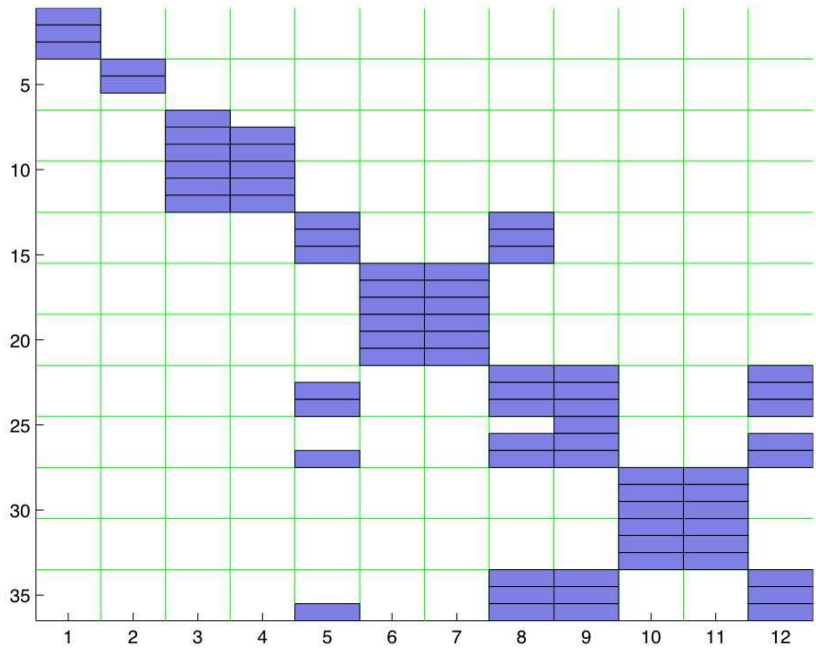

FIG 4. *0.99-confidence predction regions for tobacco "smells".*

### *3.3. Nonparametric p-values via permutation tests*

For a particular class $\theta$ let $I(1) < I(2) < \cdots < I(N_\theta)$ be the elements of $\mathcal{G}_\theta$. An elementary but useful fact is that $(X, X_{I(1)}, X_{I(2)}, \ldots, X_{I(N_\theta)})$ is exchangeable conditional on $Y = \theta$. Thus let $T_\theta(X, \mathcal{D})$ be a test statistic which is symmetric in $(X_{I(j)})_{j=1}^{N_\theta}$. We define $\mathcal{D}_i(x)$ to be the training data with $x$ in place of $X_i$.

Then the nonparametric p-value

$$\pi_\theta(X, \mathcal{D}) \ := \ \frac{\#\{i \in \mathcal{G}_\theta : T_\theta(X_i, \mathcal{D}_i(X)) \geq T_\theta(X, \mathcal{D})\} + 1}{N_\theta + 1} \quad (3.2)$$

satisfies requirement (1.2). Since $\pi_\theta \geq (N_\theta + 1)^{-1}$, this procedure is useful only if $N_\theta + 1 \geq \alpha^{-1}$. In case of $\alpha = 0.05$ this means that $N_\theta$ should be at least 19.

As for the test statistic $T_\theta(X, \mathcal{D})$, the optimal p-value in Section 2 suggests using an estimator for the weighted likelihood ratio $T_\theta^*(x)$ or a strictly increasing transformation thereof. In very high-dimensional settings this may be too ambitious, and $T_\theta(X, \mathcal{D})$ could be any test statistic quantifying the implausibility of "$Y = \theta$".

**Plug-in statistic for standard gaussian model.** For the setting of Example 2.1 and Section 3.2 one could replace the unknown parameters $w_c$, $\mu_c$ and $\Sigma$ in $T_\theta^*$ with $N_c/n$, $\widehat{\mu}_c$ and $\widehat{\Sigma}$, respectively. Note that the resulting p-values always satisfy (1.2), even if the underlying distributions $P_c$ are *not* gaussian with common covariance matrix.

**Nearest-neighbor estimation.** One could estimate $w_\theta(\cdot)$ via nearest neighbors. Suppose that $d(\cdot, \cdot)$ is some metric on $\mathcal{X}$. Let $B(x, r) := \{y \in \mathcal{X} : d(x, y) \leq r\}$, and for a fixed positive integer $k \leq n$ define

$$\widehat{r}_k(x) = \widehat{r}_k(x, \mathcal{D}) \ := \ \min\{r \geq 0 : \#\{i \leq n : X_i \in B(x, r)\} \geq k\}.$$

Further let $\widehat{P}_\theta$ denote the empirical distribution of the points $X_i$, $i \in \mathcal{G}_\theta$, i.e.

$$\widehat{P}_\theta(B) \ := \ N_\theta^{-1} \#\{i \in \mathcal{G}_\theta : X_i \in B\} \quad \text{for } B \subset \mathcal{X}.$$

Then the $k$-nearest-neighbor estimator of $w_\theta(x)$ is given by

$$\widehat{w}_\theta(x, \mathcal{D}) \ := \ \widehat{w}_\theta \widehat{P}_\theta\big(B(x, \widehat{r}_k(x))\big) \Big/ \sum_{b \in \Theta} \widehat{w}_b \widehat{P}_b\big(B(x, \widehat{r}_k(x))\big)$$

with certain estimators $\widehat{w}_b = \widehat{w}_b(\mathcal{D})$ of $w_b$. The resulting nonparametric p-value is defined with $T_\theta(x, \mathcal{D}) := -\widehat{w}_\theta(x, \mathcal{D})$. Note that in case of $\widehat{w}_b = N_b/n$, we simply end up with the ratio

$$\widehat{w}_\theta(x, \mathcal{D}) \ = \ \#\{i \in \mathcal{G}_\theta : d(X_i, x) \leq \widehat{r}_k(x)\} \Big/ \#\{i \leq n : d(X_i, x) \leq \widehat{r}_k(x)\}.$$

For simplicity, we assume $k$ to be determined by the group sizes $N_1, \ldots, N_L$ only. Of course one could define $\pi_\theta(X, \mathcal{D})$ with $k = k_\theta(X, \mathcal{D})$ nearest neighbors of $X$, as long as $k_\theta(X, \mathcal{D})$ is symmetric in the $N_\theta + 1$ feature vectors $X$ and $X_i$, $i \in \mathcal{G}_\theta$. Moreover, in applications where the different components of $X$ are measured on rather different scales, it might be reasonable to replace $d(\cdot, \cdot)$ with some data-driven metric.

**Logistic regression.** Suppose for simplicity that there are $L = 2$ classes and that $X \in \mathbb{R}^d$ contains the values of $d$ numerical or binary variables. Let $(\widehat{a}, \widehat{b}) = (\widehat{a}(\mathcal{D}), \widehat{b}(\mathcal{D}))$ be the maximum likelihood estimator for the parameter $(a, b) \in \mathbb{R} \times \mathbb{R}^d$ in the logistic model, where

$$\log \frac{w_2(x)}{1 - w_2(x)} = a + b^\top x.$$

Then possible candidates for $T_1(x, \mathcal{D})$ and $T_2(x, \mathcal{D})$ are given by

$$T_1(x, \mathcal{D}) := \widehat{a} + \widehat{b}^\top x =: -T_2(x, \mathcal{D}).$$

Extensions to multicategory logistic regression as well as the inclusion of regularization terms to deal with high-dimensional covariable vectors $X$ are possible and will be described elsewhere.

### 3.4. Asymptotic properties

Now we analyze the asymptotic behavior of the nonparametric p-values $\pi_\theta(X, \mathcal{D})$ and the corresponding empirical probabilities $\widehat{\mathcal{I}}_\alpha(b, \theta)$ and $\widehat{\mathcal{P}}(b, S)$. Throughout this section, asymptotic statements are to be understood within setting (3.1).

As in Section 2 we assume that the distributions $P_\theta$ have strictly positive densities with respect to some measure $M$ on $\mathcal{X}$. The following theorem implies that $\pi_\theta(X, \mathcal{D})$ satisfies (1.3) under certain conditions on the underlying test statistic $T_\theta(X, \mathcal{D})$. In addition the empirical probabilities $\widehat{\mathcal{I}}_\alpha(b, \theta)$ and $\widehat{\mathcal{P}}(b, S)$ turn out to be consistent estimators of $\mathcal{I}_\alpha(b, \theta \mid \mathcal{D})$ and $\mathcal{P}_\alpha(b, S \mid \mathcal{D})$, respectively.

**Theorem 3.1.** *Suppose that for fixed $\theta \in \Theta$ there exists a test statistic $T_\theta^o$ on $\mathcal{X}$ satisfying the following two requirements:*

$$T_\theta(X, \mathcal{D}) \to_p T_\theta^o(X), \tag{3.3}$$

$$\mathcal{L}(T_\theta^o(X)) \text{ is continuous.} \tag{3.4}$$

*Then*

$$\pi_\theta(X, \mathcal{D}) \to_p \pi_\theta^o(X), \tag{3.5}$$

*where*

$$\pi_\theta^o(x) := P_\theta\{z \in \mathcal{X} : T_\theta^o(z) \geq T_\theta^o(x)\}.$$

*In particular, for arbitrary fixed $\alpha \in (0, 1)$,*

$$\mathcal{R}_\alpha(\pi_\theta(\cdot, \mathcal{D})) \to_p \mathcal{R}_\alpha(\pi_\theta^o), \tag{3.6}$$

$$\left.\begin{array}{c}\mathcal{I}_\alpha(b, \theta \mid \mathcal{D}) \\ \widehat{\mathcal{I}}_\alpha(b, \theta)\end{array}\right\} \to_p \mathbb{P}(\pi_\theta^o(X) > \alpha \mid Y = b) \quad \text{for each } b \in \Theta. \tag{3.7}$$

If the limiting test statistic $T_\theta^o$ is equal to $T_\theta^*$ or some strictly increasing transformation thereof, then the nonparametric p-value $\pi_\theta(\cdot, \mathcal{D})$ is asymptotically optimal. The next two lemmata describe situations in which Condition (3.3) or (3.4) is satisfied.

**Lemma 3.2.** *Conditions (3.3) and (3.4) are satisfied in case of the plug-in rule for the homoscedastic gaussian model, provided that $\mathbb{E}(\|X\|^2) < \infty$ and $\mathcal{L}(X)$ has a Lebesgue density.*

**Lemma 3.3.** *Suppose that $(\mathcal{X}, d)$ is a separable metric space and that all densities $f_b$, $b \in \Theta$, are continuous on $\mathcal{X}$. Alternatively, suppose that $\mathcal{X} = \mathbb{R}^q$ equipped with some norm. Then Condition (3.3) is satisfied with $T_\theta^o = T_\theta^*$ in case of the k-nearest-neighbor rule with $\widehat{w}_\theta = N_\theta/n$, provided that*

$$k = k(n) \to \infty \quad \text{and} \quad k/n \to 0.$$

### 3.5. Examples

The nonparametric p-values are illustrated with two examples.

**Example 3.2.** The lower right panel in Figure 1 shows simulated training data from the model in Example 2.2, where $N_1 = N_2 = N_3 = 100$. Now we computed the corresponding prediction regions $\widehat{\mathcal{Y}}_{0.05}(x, \mathcal{D})$ based on the plug-in method for the standard gaussian model (which isn't correct here) and on the nearest-neighbor method with $k = 100$ and standard euclidean distance. Figure 5 depicts these prediction regions.

To judge the performance of the nonparametric p-values visually we chose ROC curves, where we concentrated on the plug-in method. In Figure 6 we show for each pair $(b, \theta) \in \Theta \times \Theta$ the true ROC curves of $\pi_\theta^*(\cdot)$ and $\pi_\theta(\cdot, \mathcal{D})$,

$$(0,1) \ni \alpha \mapsto \begin{cases} \mathbb{P}\big(\pi_\theta^*(X) \leq \alpha \,\big|\, Y = b\big) & \text{(magenta)}, \\ \mathbb{P}\big(\pi_\theta(X, \mathcal{D}) \leq \alpha \,\big|\, Y = b, \mathcal{D}\big) = 1 - \mathcal{I}_\alpha(b, \theta \,|\, \mathcal{D}) & \text{(blue)}, \end{cases}$$

both of which had been estimated in 40'000 Monte Carlo Simulations of $X \sim P_\theta$. In addition we show the empirical ROC curve $\alpha \mapsto 1 - \widehat{\mathcal{I}}_\alpha(b, \theta)$ (black step function). Note first that the difference between the (conditional) ROC curve of $\pi_\theta(\cdot, \mathcal{D})$ and its empirical counterpart $1 - \mathcal{I}_\alpha(b, \theta \,|\, \mathcal{D})$ is always rather small, despite the moderate group sizes $N_b = 100$. Note further that the ROC curves of $\pi_\theta(\cdot, \mathcal{D})$ and $\pi_\theta^*(\cdot)$ are also close together, despite the fact that the plug-in method uses an incorrect model. These pictures show clearly that distinguishing between classes 1 and 2 is more difficult than distinguishing between classes 2 and 3, while classes 1 and 3 are separated almost perfectly.

Of course these pictures give only partial information about the performance of the p-values. In addition one could investigate the joint distribution of the p-values via pattern probabilities; see also the next example.

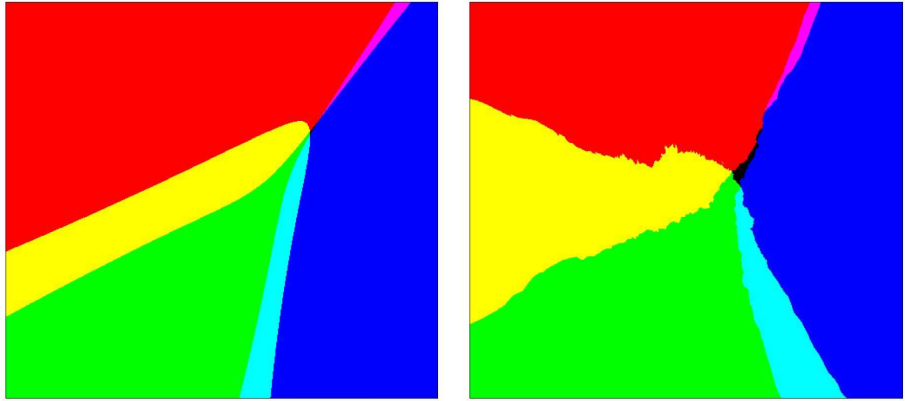

FIG 5. *Prediction regions $\widehat{\mathcal{Y}}_{0.05}(x, \mathcal{D})$ with plug-in method (left) and nearest neighbor method (right) for Example 3.2.*

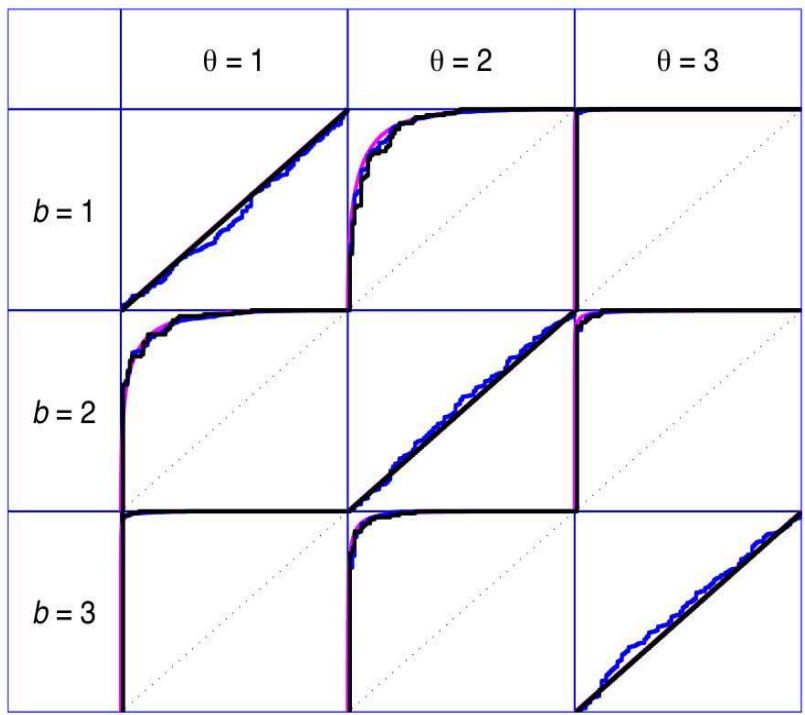

FIG 6. *ROC curves for the plug-in method applied to the data in Example 3.2.*

TABLE 1
Empirical performance of $\widehat{\mathcal{Y}}_{0.05}(\cdot,\cdot)$ and $\widehat{\mathcal{Y}}_{0.01}(\cdot,\cdot)$ in Example 3.3.

| $Y_i$ | $\widehat{\mathcal{Y}}_{0.05}(X_i,\mathcal{D}_i)$ | | | | | $\widehat{\mathcal{Y}}_{0.01}(X_i,\mathcal{D}_i)$ | | | | |
|---|---|---|---|---|---|---|---|---|---|---|
| | $\ni 1$ | $\ni 2$ | $=\{1\}$ | $=\{2\}$ | $=\{1,2\}$ | $\ni 1$ | $\ni 2$ | $=\{1\}$ | $=\{2\}$ | $=\{1,2\}$ |
| 1 | .950 | .244 | .756 | .050 | .194 | .990 | .448 | .552 | .010 | .438 |
|   | .950 | .222 | .778 | .050 | .172 | .990 | .452 | .548 | .010 | .443 |
|   | .952 | .233 | .767 | .048 | .185 | .990 | .449 | .551 | .010 | .440 |
| 2 | .396 | .950 | .050 | .604 | .346 | .743 | .991 | .009 | .257 | .734 |
|   | .356 | .950 | .050 | .644 | .307 | .698 | .991 | .009 | .302 | .689 |
|   | .406 | .950 | .050 | .594 | .356 | .773 | .992 | .008 | .227 | .766 |

**Example 3.3.** This example is from a data base on quality management at the University hospital at Lübeck. In a longterm study on mortality of patients after a certain type of heart surgery, data of more than 20'000 cases have been reported. The dependent variable is $Y \in \{1,2\}$ with $Y = 1$ and $Y = 2$ meaning that the patient survived the operation or not, respectively. For each case there were $q = 21$ numerical or binary covariables describing the patient (e.g. sex, age, various specific risk factors) plus covariables describing the circumstances of the operation (e.g. emergency or not, experience of the surgeon).

We reduced the data set by taking all $N_1 = 662$ observations with $Y = 2$ and a random subsample of $N_1 = 3N_2 = 1986$ observations with $Y = 1$. Without such a reduction, the nearest-neighbor method wouldn't work well due to the very different group sizes. Now we computed nonparametric crossvalidated p-values based on the plug-in method from the standard gaussian model, logistic regression, and the nearest-neighbor method with $k = 200$. In the latter case, we first divided each component of $X$ corresponding to a non-dichotomous variable by its sample standard deviation, because the variables are measured on very different scales. Table 1 reports the performance of $\widehat{\mathcal{Y}}_\alpha(X_i, \mathcal{D}_i)$ as a predictor of $Y_i$ for $\alpha = 5\%$ and $\alpha = 1\%$. In each cell of the table the entries correspond to the three methods mentioned above. This example shows the p-values' potential to classify a certain fraction of cases unambiguously even in situations in which overall risks of classifiers are not small which is rather typical in medical applications. Note again that the method doesn't require any knowledge of prior probabilities. Logistic regression yielded slightly better results than the other two in terms of the fraction of cases with $\widehat{\mathcal{Y}}_\alpha(X_i, \mathcal{D}_i) = \{Y_i\}$. The other two methods performed similarly.

### 3.6. Impossibility of strengthening (1.3)

Comparing (1.2) and (1.3), one might want to strengthen the latter requirement to
$$\mathbb{P}\big(\pi_\theta(X,\mathcal{D}) \le \alpha \,\big|\, Y = \theta, \mathcal{D}\big) \le \alpha \quad \text{almost surely.} \tag{3.8}$$
However, the following lemma entails that there are no reasonable p-values satisfying (3.8). Recall that we are aiming at p-values such that $\mathbb{P}\big(\pi_\theta(X,\mathcal{D}) \le \alpha \,\big|\, Y = b\big)$ is large for $b \ne \theta$.

**Lemma 3.4.** *Let $Q_1, Q_2, \ldots, Q_L$ be mutually absolutely continuous probability distributions on $\mathcal{X}$. Suppose that (3.8) is satisfied whenever $(P_1, P_2, \ldots, P_L)$ is a permutation of $(Q_1, Q_2, \ldots, Q_L)$. In that case, for arbitrary $b \in \Theta$,*

$$\mathbb{P}\big(\pi_\theta(X, \mathcal{D}) \leq \alpha \,\big|\, Y = b, \mathcal{D}\big) \;\leq\; \alpha \quad \textit{almost surely.}$$

## 4. Computational aspects

The computation of the p-values in (3.2) may be rather time-consuming, depending on the particular test statistic $T_\theta(\cdot, \mathcal{D})$. Just think about classification methods involving variable selection or tuning of artificial neural networks by means of $\mathcal{D}$. Also the nearest-neighbor method with some data-driven choice of $k$ or the metric $d(\cdot, \cdot)$ may result in tedious procedures. In order to compute $\pi_\theta(\cdot, \mathcal{D})$ as well as $\pi_\theta(X_i, \mathcal{D}_i)$ one can typically reduce the computational complexity considerably by using suitable update formulae or shortcuts.

**Naive shortcuts for the nonparametric p-values.** One might be tempted to replace $\pi_\theta(X, \mathcal{D})$ with the naive p-values

$$\pi_\theta^{\mathrm{naive}}(X, \mathcal{D}) \;:=\; \frac{\#\{i \in \mathcal{G}_\theta : T_\theta(X_i, \mathcal{D}) \geq T_\theta(X, \mathcal{D})\} + 1}{N_\theta + 1}. \tag{4.1}$$

One can easily show that the conclusions of Theorem 3.1 remain true with $\pi_\theta^{\mathrm{naive}}(\cdot, \cdot)$ in place of $\pi_\theta(\cdot, \cdot)$. However, finite sample validity in the sense of (1.2) is not satisfied in general, so we prefer the alternative shortcut described next. Note also that empirical ROC curves offered by some statistical software packages, as a complement to logistic regression or linear discriminant analysis with two classes, are often based on this shortcut.

**Valid shortcuts for the nonparametric p-values.** Often the computations as well as the program code become much simpler if we replace $T_\theta(X, \mathcal{D})$ and $T_\theta(X_i, \mathcal{D}_i(X))$ in Definition (3.2) with $T_\theta(X, \mathcal{D}(X, \theta))$ and $T_\theta(X_i, \mathcal{D}(X, \theta))$, respectively, where $\mathcal{D}(X, \theta)$ denotes the training data $\mathcal{D}$ after *adding* the "observation" $(X, \theta)$. That means, before judging whether $\theta$ is a plausible class label for a new observation $X$, we augment the training data by $(X, \theta)$ to determine the test statistic $T_\theta(\cdot, \mathcal{D}(X, \theta))$. Then we just evaluate the latter function at the $N_\theta + 1$ points $X$ and $X_i$, $i \in \mathcal{G}_\theta$, to compute

$$\pi_\theta^{\mathrm{naive}}(X, \mathcal{D}(X, \theta)) \;=\; \frac{\#\{i \in \mathcal{G}_\theta : T_\theta(X_i, \mathcal{D}(X, \theta)) \geq T_\theta(X, \mathcal{D}(X, \theta))\} + 1}{N_\theta + 1}.$$

This p-value *does* satisfy Condition (1.2), and the conclusions of Theorem 3.1 remain true as well. In this context it might be helpful if the underlying test statistics satisfy some moderate robustness properties, because $X$ may be an outlier with respect to the distribution $P_\theta$.

**Update formulae for sample means and covariances.** In connection with the typicality indices of Section 3.2 or the plug-in method for the standard gaussian model, elementary calculations reveal the following update formulae for groupwise mean vectors and sample covariance matrices: Replacing $\mathcal{D}$ with the reduced data set $\mathcal{D}_i$ for some $i \in \mathcal{G}_\theta$ has no impact on $\widehat{\mu}_b$ for $b \neq \theta$ while

$$\widehat{\Sigma} \leftarrow (n - L - 1)^{-1}\Big((n - L)\widehat{\Sigma} - (1 - N_\theta^{-1})^{-1}(X_i - \widehat{\mu}_\theta)(X_i - \widehat{\mu}_\theta)^\top\Big),$$
$$\widehat{\mu}_\theta \leftarrow (N_\theta - 1)^{-1}(N_\theta \widehat{\mu}_\theta - X_i) = \widehat{\mu}_\theta - (N_\theta - 1)^{-1}(X_i - \widehat{\mu}_\theta).$$

Replacing $\mathcal{D}$ with the modified data set $\mathcal{D}_i(X)$ for some $i \in \mathcal{G}_\theta$ results in

$$\widehat{\Sigma} \leftarrow (n - L)^{-1}\Big((n - L)\widehat{\Sigma}$$
$$+ (1 - N_\theta^{-1})\big((X - \widehat{\mu}_{\theta,i})(X - \widehat{\mu}_{\theta,i})^\top - (X_i - \widehat{\mu}_{\theta,i})(X_i - \widehat{\mu}_{\theta,i})^\top\big)\Big),$$
$$\widehat{\mu}_\theta \leftarrow \widehat{\mu}_\theta + N_\theta^{-1}(X - X_i),$$

where $\widehat{\mu}_{\theta,i} := (N_\theta - 1)^{-1}(N_\theta \widehat{\mu}_\theta - X_i)$. Finally, replacing $\mathcal{D}$ with the augmented data set $\mathcal{D}(X, \theta)$ means that

$$\widehat{\Sigma} \leftarrow (n + 1 - L)^{-1}\Big((n - L)\widehat{\Sigma} + (1 + N_\theta^{-1})^{-1}(X - \widehat{\mu}_\theta)(X - \widehat{\mu}_\theta)^\top\Big),$$
$$\widehat{\mu}_\theta \leftarrow (N_\theta + 1)^{-1}(N_\theta \widehat{\mu}_\theta + X) = \widehat{\mu}_\theta + (N_\theta + 1)^{-1}(X - \widehat{\mu}_\theta).$$

**Update formulae for the nearest-neighbor method.** For convenience we restrict our attention to the valid shortcut involving $\mathcal{D}(X, \theta)$. To compute the resulting p-values $\widehat{\pi}_\theta^{\text{naive}}(X, \mathcal{D}(X, \theta))$ quickly for arbitrary feature vectors $X \in \mathcal{X}$, it is convenient to store the $n(1 + 2L)$ numbers

$$\widehat{r}_k(X_i, \mathcal{D}), \ N_{k-1,b}(X_i, \mathcal{D}), \ N_{k,b}(X_i, \mathcal{D})$$

with $i \in \{1, \ldots, n\}$ and $b \in \Theta$, where

$$N_{\ell,b}(x, \mathcal{D}) := \#\big\{i \in \{1, \ldots, n\} : Y_i = b, d(x, X_i) \leq \widehat{r}_\ell(x, \mathcal{D})\big\}.$$

For then one can easily verify that

$$N_{k,b}(X_i, \mathcal{D}(X, \theta)) = \begin{cases} N_{k-1,b}(X_i, \mathcal{D}) + 1\{b = \theta\} & \text{if } d(X_i, X) < \widehat{r}_k(X_i, \mathcal{D}) \\ N_{k,b}(X_i, \mathcal{D}) + 1\{b = \theta\} & \text{if } d(X_i, X) = \widehat{r}_k(X_i, \mathcal{D}), \\ N_{k,b}(X_i, \mathcal{D}) & \text{if } d(X_i, X) > \widehat{r}_k(X_i, \mathcal{D}). \end{cases}$$

Hence classifying a new feature vector $X$ requires only $O(n)$ steps for determining the $1 + L^2$ numbers $\widehat{r}_k(X, \mathcal{D}(X, \theta))$ and $N_b(X, \mathcal{D}(X, \theta))$ and the $nL^2$ numbers $N_b(X_i, \mathcal{D}(X, \theta))$, where $1 \leq i \leq n$ and $b, \theta \in \Theta$.

Computing the crossvalidated p-values with the valid shortcut is particularly easy, because replacing one training observation $(X_i, Y_i)$ with $(X_i, \theta)$ does not affect the radii $\widehat{r}_k(x, \mathcal{D})$.

In case of data-driven choice of $k$ or $d(\cdot, \cdot)$, the preceding formulae are no longer applicable. Then the valid shortcut is particularly useful to reduce the computational complexity.

## 5. Likelihood ratios and local identifiability

In previous sections we assumed that the distribution of likelihood ratios such as $w_\theta(X)$ or $T_\theta^*(X)$ is continuous. This property is related to a property which we call 'local identifiability', a strengthening of the well-known notion of identifiability for finite mixtures. Throughout this section we assume that the distributions $P_1, P_2, \ldots, P_L$ belong to a given model $(Q_\xi)_{\xi \in \Xi}$ of probability distributions $Q_\xi$ with densities $g_\xi > 0$ with respect to some measure $M$ on $\mathcal{X}$.

**Identifiability.** Let us first recall Yakowitz and Spragins' [13] definition of identifiability for finite mixtures. The family $(Q_\xi)_{\xi \in \Xi}$ is called *identifiable*, if the following condition is satisfied: For arbitrary $m \in \mathbb{N}$ let $\xi(1), \ldots, \xi(m)$ be pairwise different parameters in $\Xi$ and let $\lambda_1, \ldots, \lambda_m > 0$. If $\xi'(1), \ldots, \xi'(m) \in \Xi$ and $\lambda_1', \ldots, \lambda_m' \geq 0$ such that

$$\sum_{j=1}^m \lambda_j Q_{\xi(j)} = \sum_{j=1}^m \lambda_j' Q_{\xi'(j)},$$

then there exists a permutation $\sigma$ of $\{1, 2, \ldots, m\}$ such that $\xi'(i) = \xi(\sigma(i))$ and $\lambda_i' = \lambda_{\sigma(i)}$ for $i = 1, 2, \ldots, m$.

Evidently the family $(Q_\xi)_{\xi \in \Xi}$ is identifiable if the density functions $g_\xi$, $\xi \in \Xi$, are linearly independent as elements of $L^1(M)$, and the converse statement is also true [13].

A standard example of an identifiable family is the set of all nondegenerate gaussian distributions on $\mathbb{R}^q$; see [13]. Holzmann et al. [6] provide a rather comprehensive list of identifiable classes of multivariate distributions. In particular, they verify identifiability of families of elliptically symmetric distributions on $\mathcal{X} = \mathbb{R}^q$ with Lebesgue densities of the form

$$g_\xi(x) = \det(\Sigma)^{-1/2} h_q\big((x-\mu)^\top \Sigma^{-1}(x-\mu); \zeta\big). \tag{5.1}$$

Here the parameter $\xi = (\mu, \Sigma, \zeta)$ consists of an arbitrary location parameter $\mu \in \mathbb{R}^q$, an arbitrary symmetric and positive definite scatter matrix $\Sigma \in \mathbb{R}^{q \times q}$ and an additional shape parameter $\zeta$ which may also vary in the mixture. For each shape parameter $\zeta$, the 'density generator' $h_q(\cdot; \zeta)$ is a nonnegative function on $[0, \infty)$ such that $\int_{\mathcal{X}} h_q(\|x\|^2; \zeta) \, dx = 1$. One particular example are the multivariate $t$–distributions with

$$h_q(u; \zeta) = \frac{\Gamma((\zeta+q)/2)}{\pi^{q/2} \Gamma(\zeta/2)} (1+u)^{-(\zeta+q)/2}$$

for $\zeta > 0$. We mention that the subsequent arguments apply to most of the elliptically symmetric families discussed by Holzmann et al. [6]. Peel et al. [9] discuss classification for directional data and our method can be extended to distributions with non-euclidean domain, combining the arguments below with methods in Holzmann et al. [5]. As prominent examples we mention the von Mises family for directional data and the Kent family for spherical data.

**Continuity of likelihood ratios.** Suppose that $P_\theta = Q_{\xi(\theta)}$ with parameters $\xi(1), \ldots, \xi(L)$ in $\Xi$ which are not all identical. Then one can easily verify that continuity of $\mathcal{L}(w_\theta(X))$ or $\mathcal{L}(T_\theta^*(X))$ follows from the following condition:
The family $(Q_\xi)_{\xi \in \Xi}$ is called *locally identifiable*, if for arbitrary $m \in \mathbb{N}$, pairwise different parameters $\xi(1), \ldots, \xi(m) \in \Xi$ and numbers $\beta_1, \ldots, \beta_m \in \mathbb{R}$,

$$M\left\{x \in \mathcal{X} : \sum_{j=1}^m \beta_j g_{\xi(j)}(x) = 0\right\} > 0 \quad \text{implies that} \quad \beta_1 = \beta_2 = \cdots = \beta_m = 0.$$

Local identifiability entails the following conclusion: Suppose that $Q$ is equal to $\sum_{j=1}^m \lambda_j Q_{\xi(j)}$ for some number $m \in \mathbb{N}$, pairwise different parameters $\xi(1), \ldots, \xi(m)$ in $\Xi$ and nonnegative numbers $\lambda_1, \ldots, \lambda_m$. Then one can determine the ingredients $m$, $\xi(1), \ldots, \xi(m)$ and $\lambda_1, \ldots, \lambda_m$ from the restriction of $Q$ to any fixed measurable set $B_o \subset \mathcal{X}$ with $M(B_o) > 0$. The following theorem provides a sufficient criterion for local identifiability which is easily verified in many standard examples.

**Theorem 5.1.** *Let $M$ be Lebesgue measure on $\mathcal{X} = \mathcal{X}_1 \times \mathcal{X}_2 \times \cdots \times \mathcal{X}_q$ with open intervals $\mathcal{X}_k \subset \mathbb{R}$. Suppose that the following two conditions are satisfied:*
*(i) $(Q_\xi)_{\xi \in \Xi}$ is identifiable;*
*(ii) for arbitrary $\xi \in \Xi$, $k \in \{1, 2, \ldots, q\}$ and $x_i \in \mathcal{X}_i$, $i \neq k$, the function*

$$t \mapsto g_\xi(x_1, \ldots, x_{k-1}, t, x_{k+1}, \ldots, x_q)$$

*may be extended to a holomorphic function on some open subset of $\mathbb{C}$ containing $\mathcal{X}_k$.*
*Then the family $(Q_\xi)_{\xi \in \Xi}$ is locally identifiable.*

One can easily verify that Condition (ii) of Theorem 5.1 is satisfied by the densities $g_\xi$ in (5.1), if the density generators $h_q(\cdot; \zeta)$ may be extended to holomorphic functions on some open subset of $\mathbb{C}$ containing $[0, \infty)$. Hence, for instance, the family of all multivariate $t$–distributions is locally identifiable.

## 6. Proofs

*Proof of Theorem 3.1.* Since the distributions $P_1, \ldots, P_L$ are mutually absolutely continuous, Condition (3.3) entails that

$$\rho(\epsilon, N_1, \ldots, N_L)$$
$$:= \max_{a,b \in \Theta;\, i=1,\ldots,n} \int \mathbb{P}\bigl(|T_\theta(x, \mathcal{D}_i(z)) - T_\theta^o(x)| \geq \epsilon\bigr)\, P_a(dx) P_b(dz)$$

tends to zero for any fixed $\epsilon > 0$.

It follows from the elementary inequality

$$\bigl|1\{r \geq s\} - 1\{r_o \geq s_o\}\bigr| \leq 1\{|r - r_o| \geq \epsilon\} + 1\{|s - s_o| \geq \epsilon\} + 1\{|r_o - s_o| < 2\epsilon\}$$

for real numbers $r, r_o, s, s_o$ that

$$\begin{aligned}
\pi_\theta(X,\mathcal{D}) &= (N_\theta+1)^{-1}\left(1+\sum_{i\in\mathcal{G}_\theta}1\{T_\theta(X_i,\mathcal{D}_i(X))\geq T_\theta(X,\mathcal{D})\}\right) \\
&= N_\theta^{-1}\sum_{i\in\mathcal{G}_\theta}1\{T_\theta(X_i,\mathcal{D}_i(X))\geq T_\theta(X,\mathcal{D})\} + R_1 \\
&= N_\theta^{-1}\sum_{i\in\mathcal{G}_\theta}1\{T_\theta^o(X_i)\geq T_\theta^o(X)\} + R_1 + R_2(\epsilon),
\end{aligned}$$

where $|R_1|\leq (N_\theta+1)^{-1}$ and

$$\begin{aligned}
|R_2(\epsilon)| &\leq N_\theta^{-1}\#\Big\{i\in\mathcal{G}_\theta : \big|T_\theta(X_i,\mathcal{D}_i(X))-T_\theta^o(X_i)\big|\geq \epsilon\Big\} \\
&\quad + 1\Big\{\big|T_\theta(X,\mathcal{D})-T_X^o(X)\big|\geq \epsilon\Big\} \\
&\quad + N_\theta^{-1}\#\Big\{i\in\mathcal{G}_\theta : \big|T_\theta^o(X_i)-T_\theta^o(X)\big|<2\epsilon\Big\}.
\end{aligned}$$

Hence $\mathbb{E}\,|R_2(\epsilon)|\leq 2\rho(\epsilon,N_1,\ldots,N_L)+\omega(2\epsilon)\to\omega(2\epsilon)$, where

$$\omega(\delta) := \sup_{r\in\mathbb{R}} P_\theta\big\{z\in\mathcal{X} : |T_\theta^o(z)-r|<\delta\big\}\ \downarrow\ 0\quad(\delta\downarrow 0)$$

by virtue of Condition (3.4). These considerations show that

$$\pi_\theta(X,\mathcal{D}) = \widehat{F}_\theta(T_\theta^o(X)) + o_p(1) = F_\theta(T_\theta^o(X)) + o_p(1),$$

where

$$\begin{aligned}
F_\theta(r) &:= P_\theta\{z\in\mathcal{X} : T_\theta^o(z)\geq r\}, \\
\widehat{F}_\theta(r) &:= \widehat{P}_\theta\{z\in\mathcal{X} : T_\theta^o(z)\geq r\}.
\end{aligned}$$

Here we utilized the well-known fact [11] that $\|\widehat{F}_\theta - F_\theta\|_\infty = o_p(1)$. Since $\pi_\theta^o(X) = F_\theta(T_\theta^o(X))$, this entails Conclusion (3.5).

As to the remaining assertions (3.6–3.7), note first that (3.5) implies that

$$\begin{aligned}
&\tau(\epsilon, N_1,\ldots,N_L) \\
&:= \max_{a,b\in\Theta;\,i=1,\ldots,n}\int \mathbb{P}\big(|\pi_\theta(x,\mathcal{D}_i(z))-\pi_\theta^o(x)|\geq\epsilon\big)\,P_a(dx)P_b(dz)
\end{aligned}$$

tends to zero for any fixed $\epsilon>0$, again a consequence of mutual absolute continuity of $P_1,\ldots,P_L$. Similarly as in the proof of (3.5) one can verify that

$$\begin{aligned}
\mathcal{I}_\alpha(b,\theta\,|\,\mathcal{D}) = \mathbb{P}(\pi_\theta(X,\mathcal{D})>\alpha\,|\,Y=b,\mathcal{D}) &= G_{b,\theta}(\alpha)+R(\epsilon), \\
\widehat{\mathcal{I}}_\alpha(b,\theta) = N_b^{-1}\sum_{i\in\mathcal{G}_b}1\{\pi_\theta(X_i,\mathcal{D}_i)>\alpha\} &= \widehat{G}_{b,\theta}(\alpha)+\widehat{R}(\epsilon) \\
&= G_{b,\theta}(\alpha)+\widehat{R}(\epsilon)+o_p(1),
\end{aligned}$$

with $G_{b,\theta}(u) := P_b\{z \in \mathcal{X} : \pi_\theta^o(z) > u\}$ and $\widehat{G}_{b,\theta}(u) := \widehat{P}_b\{z \in \mathcal{X} : \pi_\theta^o(z) > u\}$, while

$$\begin{aligned}
\mathbb{E}\,|R(\epsilon)| &\le \tau(\epsilon, N_1, \ldots, N_L) + \mathbb{P}\big(|\pi_\theta^o(X) - \alpha| < \epsilon \,\big|\, Y = b\big) \\
&\to \mathbb{P}\big(|\pi_\theta^o(X) - \alpha| < \epsilon \,\big|\, Y = b\big), \\
\mathbb{E}\,|\widehat{R}(\epsilon)| &\le \tau(\epsilon, N_1, \ldots, N_{b-1}, N_b - 1, N_{b+1}, \ldots, N_L) \\
&\quad + \mathbb{P}\big(|\pi_\theta^o(X) - \alpha| < \epsilon \,\big|\, Y = b\big) \\
&\to \mathbb{P}\big(|\pi_\theta^o(X) - \alpha| < \epsilon \,\big|\, Y = b\big).
\end{aligned}$$

Since the latter probability tends to zero as $\epsilon \downarrow 0$, we obtain Claim (3.7).

This implies Claim (3.6), because

$$\begin{aligned}
\mathcal{R}_\alpha(\pi_\theta(\cdot, \mathcal{D})) &= \sum_{b \in \Theta} w_b \mathcal{I}_\alpha(b, \theta \,|\, \mathcal{D}) \\
&\to_p \sum_{b \in \Theta} w_b \,\mathbb{P}(\pi_\theta^o(X) > \alpha \,|\, Y = b) \;=\; \mathcal{R}_\alpha(\pi_\theta^o). \qquad \square
\end{aligned}$$

*Proof of Lemma 3.2.* It is a simple consequence of the weak law of large numbers that $\widehat{\mu}_b \to_p \mu_b := \mathbb{E}(X \,|\, Y = b)$ and $\widehat{\Sigma} \to_p \Sigma := \sum_{b=1}^L w_b \operatorname{Var}(X \,|\, Y = b)$. Now one can easily show that (3.3) is satisfied with $T_\theta^o$ defined as in (2.1). The results from Section 5 entail that $\operatorname{Leb}_q\{x \in \mathbb{R}^q : T_\theta^o(x) = c\} = 0$ for any $c > 0$, so that (3.4) is satisfied as well. $\square$

*Proof of Lemma 3.3.* The assumptions imply the existence of a Borel set $\mathcal{X}_o \subset \mathcal{X}$ with $\mathbb{P}(X \in \mathcal{X}_o) = 1$ such that the following additional requirements are satisfied:

$$\mathbb{P}(X \in B(x, r)) > 0 \quad \text{for all } x \in \mathcal{X}_o, r > 0, \tag{6.1}$$

$$\lim_{r \downarrow 0} \frac{P_b(B(x,r))}{P_\theta(B(x,r))} = \frac{f_b}{f_\theta}(x) \quad \text{for all } \theta, b \in \Theta, x \in \mathcal{X}_o. \tag{6.2}$$

In case of continuous densities $f_1, f_2, \ldots, f_L > 0$ on a separable metric space $(\mathcal{X}, d)$, this is easily verified with $\mathcal{X}_o$ being the support of $\mathcal{L}(X)$, i.e. the smallest closed set such that $\mathbb{P}(X \in \mathcal{X}_o) = 1$. In case of $\mathcal{X} = \mathbb{R}^q$ and $d(x, y) = \|x - y\|$, existence of such a set $\mathcal{X}_o$ is a known result from geometric measure theory; cf. Federer [2, Theorem 2.9.8].

In view of (6.1–6.2), it suffices to show that for arbitrary fixed $x \in \mathcal{X}_o$ and $b \in \Theta$,

$$\widehat{r}_{k(n)}(x) \to_p 0 \quad \text{and} \quad \frac{\widehat{P}_b\big(B(x, \widehat{r}_{k(n)}(x))\big)}{P_b\big(B(x, \widehat{r}_{k(n)}(x))\big)} \to_p 1. \tag{6.3}$$

To this end, note first that the random numbers $N(x, r) := \#\{i : d(X_i, x) < r\}$

satisfy

$$\begin{aligned}
\mathbb{E}\, N(x,r) &= \sum_{\theta \in \Theta} N_\theta P_\theta\{z : d(z,x) < r\} \\
&= n\bigl(\mathbb{P}(d(X,x) < r) + o(1)\bigr) \quad \text{uniformly in } r \geq 0, \quad (6.4)\\
\mathrm{Var}(N(x,r)) &= \sum_{\theta \in \Theta} N_\theta P_\theta\{z : d(z,x) < r\}\bigl(1 - P_\theta\{z : d(z,x) < r\}\bigr)\\
&\leq \min\bigl\{\mathbb{E}\, N(x,r), n/4\bigr\}. \quad (6.5)
\end{aligned}$$

If we define $r_n := \max\{r \geq 0 : \mathbb{E}\, N(x,r) \leq k(n)/2\}$, then

$$\begin{aligned}
\mathbb{P}\bigl(\widehat{r}_{k(n)}(x) < r_n\bigr) &= \mathbb{P}\bigl(N(x,r_n) \geq k(n)\bigr)\\
&\leq \mathbb{P}\bigl(N(x,r_n) - \mathbb{E}\, N(x,r_n) \geq k(n)/2\bigr)\\
&\leq \mathbb{E}\, N(x,r_n)/(k(n)/2)^2\\
&\leq 2/k(n) \;\to\; 0
\end{aligned}$$

by Tshebyshev's inequality and (6.5). On the other hand, for any fixed $\epsilon > 0$,

$$\begin{aligned}
\mathbb{P}\bigl(\widehat{r}_{k(n)}(x) \geq \epsilon\bigr) &= \mathbb{P}\bigl(N(x,\epsilon) < k(n)\bigr)\\
&= \mathbb{P}\Bigl(N(x,\epsilon) - \mathbb{E}\, N(x,\epsilon) \leq n\bigl(o(1) - \mathbb{P}(d(X,x) < \epsilon)\bigr)\Bigr)\\
&= O(1/n)
\end{aligned}$$

according to (6.4) and (6.1). These considerations show that $\widehat{r}_{k(n)}(x) \to_p 0$, but $\widehat{r}_{k(n)}(x) \geq r_n$ with asymptotic probability one. Now we utilize that the process

$$r \;\mapsto\; \frac{\widehat{P}_b(B(x,r))}{P_b(B(x,r))} - 1$$

is a zero mean reverse martingale on $\{r \geq 0 : \mathbb{P}(d(X,x) \leq r) > 0\} \supset (0,\infty)$, so that Doob's inequality entails that

$$\mathbb{E} \sup_{r \geq r_n} \left| \frac{\widehat{P}_b(B(x,r))}{P_b(B(x,r))} - 1 \right|^2 \;\leq\; \frac{4}{N_b P_b(B(x,r_n))} \;=\; O(k(n)^{-1});$$

see Shorack and Wellner [11, Sections 3.6 and A.10-11]. The latter considerations imply the second part of (6.3). □

*Proof of Theorem 5.1.* The proof is by contradiction. To this end suppose that there are $m \geq 2$ pairwise different parameters $\xi(1), \xi(2), \ldots, \xi(m) \in \Xi$ and nonzero real numbers $\beta_1, \beta_2, \ldots, \beta_m$ such that $h := \sum_{i=1}^m \beta_i g_{\xi(i)}$ satisfies

$$\mathrm{Leb}_q(W) > 0 \quad \text{with} \quad W := \{x \in \mathcal{X} : h(x) = 0\}.$$

In case of $q = 1$, this entails that $W \subset \mathcal{X} = \mathcal{X}_1$ contains an accumultation point within $\mathcal{X}_1$, and the identity theorem for analytic functions yields that $h = 0$ on $\mathcal{X}$. But this would be a contradiction to $(Q_\xi)_{\xi \in \Xi}$ being identifiable.

In case of $q > 1$, by Fubini's theorem,

$$\text{Leb}_q(W) = \int_{\mathcal{X}_1 \times \cdots \times \mathcal{X}_{q-1}} \text{Leb}_1\{t : (x', t) \in W\} \, \text{Leb}_{q-1}(dx') > 0,$$

whence $\text{Leb}_1\{t : (x', t) \in W\} > 0$ for all $x'$ in a measurable set $W' \subset \mathcal{X}_1 \times \cdots \times \mathcal{X}_{q-1}$ such that $\text{Leb}_{q-1}(W') > 0$. Hence the identity theorem for analytic functions, applied to $t \mapsto h(x', t)$ implies that

$$W' \times \mathcal{X}_q \subset W.$$

Since $\text{Leb}_{q-1}(W') > 0$, we may proceed inductively, considering for $k = q-1, q-2, \ldots, 1$ the functions $t \mapsto h(x'', t, x_{k+1}, \ldots, x_q)$ on $\mathcal{X}_k$. Eventually we obtain $W = \mathcal{X}$, but this would be a contradiction to $(Q_\xi)_{\xi \in \Xi}$ being identifiable. $\square$

*Proof of Lemma 3.4.* For any permutation $\sigma$ of $(1, 2, \ldots, L)$ let $\mathbb{P}_\sigma(\cdot)$ and $\mathcal{L}_\sigma(\cdot)$ denote probabilities and distributions in case of $P_b = Q_{\sigma(b)}$ for $b = 1, 2, \ldots, L$. By assumption (3.8), for any such $\sigma$ there is a set $\mathcal{A}_\sigma$ of potential training data sets $\mathcal{D}$ such that $\mathbb{P}_\sigma(\mathcal{D} \in \mathcal{A}_\sigma) = 1$ and

$$\int 1\{\pi_\theta(x, \mathcal{D}) \leq \alpha\} \, Q_{\sigma(\theta)}(dx) \leq \alpha \quad \text{whenever } \mathcal{D} \in \mathcal{A}_\sigma.$$

Since the $L!$ distributions $\mathcal{L}_\sigma(\mathcal{D})$ are mutually absolutely continuous, the intersection $\mathcal{A} := \bigcap_\sigma \mathcal{A}_\sigma$ satisfies $\mathbb{P}_\sigma(\mathcal{D} \in \mathcal{A}) = 1$ for any permutation $\sigma$. But then

$$\int 1\{\pi_\theta(x, \mathcal{D}) \leq \alpha\} \, Q_b(dx) \leq \alpha \quad \text{for all } b \in \Theta, \mathcal{D} \in \mathcal{A}.$$

This implies that $\mathbb{P}(\pi_\theta(X, \mathcal{D}) \leq \alpha \,|\, Y = b, \mathcal{D}) \leq \alpha$ almost surely for all $b \in \Theta$, provided that $(P_1, \ldots, P_L)$ is a permutation of $(Q_1, \ldots, Q_L)$. $\square$

## Acknowledgements

We are indebted to Wolf Münchmeyer and his colleagues from *Airsense* (Schwerin) and C. Bürk (Lübeck) for fruitful conversations about classification and the data in examples 3.1 and 3.3. We are also grateful to Jerome Friedman, Trevor Hastie and Robert Tibshirani for stimulating discussions, and to Larry Wasserman for constructive comments. Lars Hömke kindly supported us in implementing some of the p-values.

## References

[1] EHM, W., E. MAMMEN and D.W. MÜLLER (1995). Power robustification of approximately linear tests. *J. Amer. Statist. Assoc.* **90**, 1025–1033. MR1354019

[2] Federer, H. (1969). *Geometric Measure Theory.* Springer, Berlin Heidelberg. MR0257325
[3] Fisher, R.A. (1936). The use of multiple measurements in taxonomic problems. *Ann. Eugenics* **7**, 179–184.
[4] Fraley, C. and A.E. Raftery (2002). Model-based clustering, discriminant analysis and density estimation. *J. Amer. Statist. Assoc.* **97**, 611–631. MR1951635
[5] Holzmann, H., A. Munk and B. Stratmann (2004). Identifiability of finite mixtures - with applications to circular distributions. *Sankhya* **66**, 440–450. MR2108200
[6] Holzmann, H., A. Munk and T. Gneiting (2006). Identifiability of finite mixtures of elliptical distributions. *Scand. J. Statist.* **33**, 753-763. MR2300914
[7] McLachlan, G.J. (1992). *Discriminant Analysis and Statistical Pattern Recognition.* Wiley, New York. MR1190469
[8] Peel, D. and G.J. McLachlan (2000). Robust mixture modeling using the $t$-distribution. *Statist. Computing* **10**, 339–348.
[9] Peel, D., W.J. Whitten and G.J. McLachlan (2001). Fitting mixtures of Kent distributions to aid in joint set identification. *J. Amer. Statist. Assoc.* **96**, 56–63. MR1973782
[10] Ripley, B.D. (1996). *Pattern Recognition and Neural Networks.* Cambridge University Press, Cambridge, UK. MR1438788
[11] Shorack, G.R. and J.A. Wellner (1986). *Empirical Processes with Applications to Statistics.* Wiley, New York. MR0838963
[12] Stone, C.J. (1977). Consistent nonparametric regression. *Ann. Statist.* **5**, 595–645. MR0443204
[13] Yakowitz, S.J. and J.D. Spragins (1968). On the identifiability of finite mixtures. *Ann. Math. Statist.* **39**, 209–214. MR0224204